\documentclass[a4paper,11pt]{article}
\usepackage{amsfonts}
\usepackage{amstext}
\usepackage{amsthm}
\usepackage{amsmath}
\usepackage{amssymb}
\usepackage{graphicx}
\usepackage{subfig}
\usepackage{latexsym}
\usepackage[latin1]{inputenc}
\newcommand{\R}{\Bbb{R}}

\newtheorem{teor}{Theorem}[section]
\newtheorem{propo}{Proposition}[section]

\newtheorem{lema}{Lemma}[section]

\newtheorem{rem}{Remark}[section]

\newcommand{\n}{\noindent}

\newcommand {\fim}{\rule{0.5em}{0.5em}}

\begin{document}

\title{On positive solutions of fractional\\
Lane-Emden systems with weights
}

\author{\textbf{Edir Junior Ferreira Leite \footnote{\textit{E-mail addresses}:
edirjrleite@ufv.br (E.J.F. Leite)}}\\ {\small\it Departamento de Matem\'{a}tica,
Universidade Federal de Vi\c{c}osa,}\\ {\small\it CCE, 36570-000, Vi\c{c}osa, MG, Brazil}
}
\date{}{

\maketitle

\markboth{abstract}{abstract}
\addcontentsline{toc}{chapter}{abstract}

\hrule \vspace{0,2cm}

\n {\bf Abstract}

In this paper we discuss the existence and regularity of solutions of the following system on a smooth bounded domain $\Omega$ in $\R^N$:
\[
\left\{
\begin{array}{llll}
\mathcal{L} u=\dfrac{v^p}{\vert x\vert^\alpha} & {\rm in} \ \ \Omega\\\\
\mathcal{L} v = \dfrac{u^q}{\vert x\vert^\beta} & {\rm in} \ \ \Omega\\\\
u= v=0 & {\rm on} \ \ \Sigma
\end{array}
\right.
\]
where $\alpha,\beta<N$ and $\mathcal{L}$ refer to any of the two types of operators $\mathcal{A}^s$ or $(-\Delta)^s$, $0 < s < 1$, 
\begin{itemize}
      \item[$\bullet$] $\Sigma = \partial\Omega$ for the spectral fractional Laplace operator $\mathcal{A}^s$,
      \item[$\bullet$] $\Sigma = \mathbb{R}^N\setminus\Omega$ for the restricted fractional Laplace operator $(-\Delta)^s$.
\end{itemize}
We find the existence of a critical hyperbole in the $(p, q)$ plane (depending on $\alpha,\beta$ and $N$) below which there exists nontrivial solutions. For such solutions, we prove an $L^\infty$ estimate of Brezis-Kato type and derive the regularity property of the solutions.
\vspace{0.5cm}
\hrule\vspace{0.2cm}

\section{Introduction and main result}

This work is devoted to the study of existence and regularity of solutions for nonlocal elliptic systems on bounded domains which will be described henceforth.

The fractional Laplace operator (or fractional Laplacian) of order $2s$, with $0 < s < 1$, denoted by $(-\Delta)^{s}$, is defined as

\[
(-\Delta)^{s}u(x) = C(N,s)\, {\rm P.V.}\int\limits_{\R^{N}}\frac{u(x)-u(y)}{\vert x-y\vert^{N+2s}}\; dy\, ,
\]

\n for all $x \in \R^{N}$, where P.V. denotes the principal value of the integral and

\[
C(N,s) = \left(\int\limits_{\R^{N}}\frac{1-\cos(\zeta_{1})}{\vert\zeta\vert^{N+2s}}\; d\zeta\right)^{-1}
\]

\n with $\zeta = (\zeta_1, \ldots, \zeta_n) \in \R^N$.

We remark that $(-\Delta)^{s}$ is a nonlocal operator on functions compactly supported in $\R^N$, i.e., to check whether the equation holds at a point, information about the values of the function far from that point is needed. 

Factional Laplace operators arise naturally in several different areas such as Probability, Finance, Physics, Chemistry and Ecology, see \cite{A, bucur}.

A closely related operator but different from $(-\Delta)^{s}$, the spectral fractional Laplace operator $\mathcal{A}^{s}$, is defined in terms of the Dirichlet spectra of the Laplace operator on $\Omega$, see \cite{compare, rafaella}. Roughly, if $(\varphi_k)$ denotes a $L^2$-orthonormal basis of eigenfunctions corresponding to eigenvalues $(\lambda_k)$ of the Laplace operator with zero Dirichlet boundary values on $\partial \Omega$, then the operator $\mathcal{A}^s$ is defined as $\mathcal{A}^{s} u = \sum_{k=1}^\infty c_k \lambda_k^s \varphi_k$, where $c_k$, $k \geq 1$, are the coefficients of the expansion $u = \sum_{k=1}^\infty c_k \varphi_k$.

An interesting interplay between the two operators occur in case of periodic solutions, or when the domain is the torus, where they coincide, see \cite{torus}. However in the case general the two operators produce very different behaviors of solutions, even when one focuses only on stable solutions, see e.g. Subsection 1.7 in \cite{serena5}.

We here are interested in studying the following problem

\begin{equation}\label{1}
\left\{
\begin{array}{llll}
\mathcal{L} u=\dfrac{v^p}{\vert x\vert^\alpha} & {\rm in} \ \ \Omega\\\\
\mathcal{L} v = \dfrac{u^q}{\vert x\vert^\beta} & {\rm in} \ \ \Omega\\\\
u= v=0 & {\rm on} \ \ \Sigma
\end{array}
\right.
\end{equation}
where $\alpha,\beta<N$ and $\mathcal{L}$ refer to any of the two types of operators $\mathcal{A}^s$ or $(-\Delta)^s$, $0 < s < 1$, 
\begin{itemize}
      \item[$\bullet$] $\Sigma = \partial\Omega$ for the spectral fractional Laplace operator $\mathcal{A}^s$,
      \item[$\bullet$] $\Sigma = \mathbb{R}^N\setminus\Omega$ for the restricted fractional Laplace operator $(-\Delta)^s$.
\end{itemize}

Our main concern in this paper is to look at the role played by the two weights when dealing with existence of solutions. We find the existence of a critical hyperbole, given by,
\begin{equation}\label{1.2}
\frac{N-\alpha}{p+1}+\frac{N-\beta}{q+1}=N-2s.
\end{equation}
Below this hyperbole we find existence of nontrivial solutions. Remark that when the two weights are not present, that is, for $\alpha=\beta=0$, we recover the critical hyperbole for elliptic systems without weights that was found independently in \cite{EM1} for $\mathcal{L}=(-\Delta)^s$ and in \cite{choi, E1} for $\mathcal{L}=\mathcal{A}^s$ (see also \cite{E2, E3, E4}). Also remark that the hyperbole (\ref{1.2}) is monotone with respect to $\alpha$ and $\beta$. The stronger the weights the smaller the hyperbole. Note also that the curve $(p,q)$ given by the hyperbole $pq = 1$ splits the behavior of (\ref{1}) into sublinear and superlinear one.

The critical hyperbole (\ref{1.2}) has vertical asymptote $p=\frac{2s-\alpha}{N-2s}$ and horizontal asymptote $q=\frac{2s-\beta}{N-2s}$. In addition, the points of intersection between the critical hyperbole (\ref{1.2}) and the hyperbole $pq = 1$ occur when $q=\frac{-(\alpha+\beta-4s)\pm\vert\alpha-\beta\vert}{2\alpha-4s}$. Thus, if $\alpha,\beta<2s$ we have such hyperboles has no point of intersection in the first quadrant (Figure 1 (a)). However, if $2s<\alpha<N$ we have such hyperboles admits an intersection point in the first quadrant (Figure 1 (b)). Similarly, if $2s<\beta<N$.

\begin{figure}[!htb]
\centering
\subfloat[Case $N>4s$ and $\alpha,\beta<2s$]{
\includegraphics[height=6.9cm]{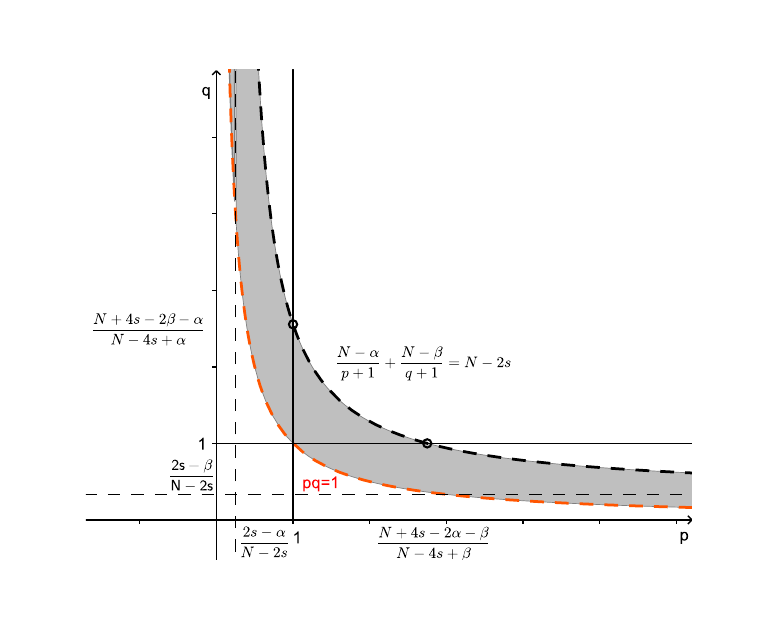}
\label{figdroopy}
}
\quad 
\subfloat[Case $2s<N<4s$, $\beta<2s$ and $2s<\alpha<N$]{
\includegraphics[height=6.9cm]{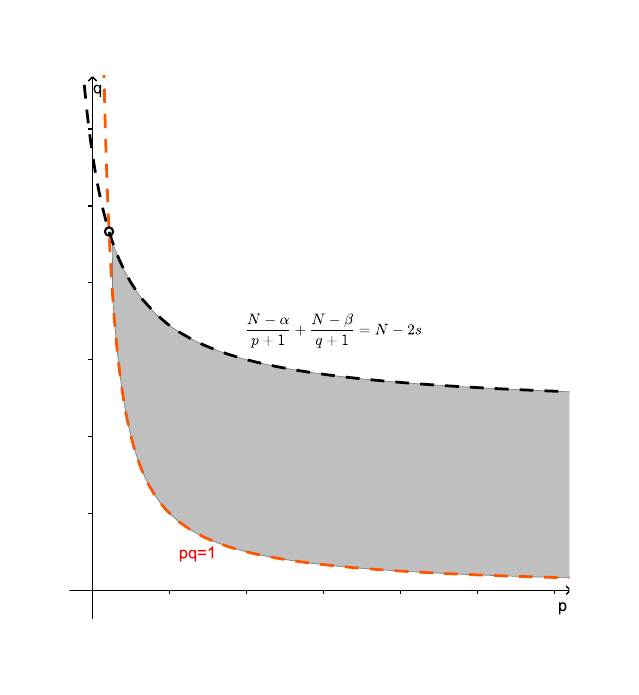}
\label{figsnoop}
}
\caption{The existence range of couples $(p,q)$}
\label{fig01}
\end{figure}

The ideas involved in our proofs base on variational methods. In particular, we use a variational argument (a linking theorem). 

By a solution of the system (\ref{1}), we mean a couple $(u,v)\in \Theta^{t}(\Omega) \times \Theta^{2s-t}(\Omega)$, satisfying
\[
\left\{
\begin{array}{llll}
\int\limits_{\Omega}\mathcal{L}^{1/2}u\mathcal{L}^{1/2}\phi dx=\int\limits_{\Omega}\dfrac{v^p}{\vert x\vert^\alpha}\phi dx & {\rm } \ \ \forall\phi\in \Theta^{2s-t}(\Omega)
\\\\
\int\limits_{\Omega}\mathcal{L}^{1/2}\psi\mathcal{L}^{1/2}v dx=\int\limits_{\Omega}\dfrac{u^q}{\vert x\vert^\beta}\psi dx & {\rm } \ \ \forall\psi\in \Theta^{t}(\Omega)
\end{array}.
\right.
\]

Our main result is

\begin{teor}\label{teo1} Let us assume that $p,q,\alpha,\beta$ verify
\begin{equation}\label{hip}
\frac{N-\alpha}{p+1}+\frac{N-\beta}{q+1}>N-2s,
\end{equation}
\begin{equation}\label{1.4}
1>\frac{1}{p+1}+\frac{1}{q+1}
\end{equation}
and let $\frac{N}{2}<t<2s$ be if $2s<N<4s$ and let $0 < t < 2s$ be such that
\[
q+1<\frac{2(N-\beta)}{N-2t} \text{ and }p+1<\frac{2(N-\alpha)}{N-(4s-2t)}\text{ if }N\geq 4s.
\]
Then, there exists infinitely many solutions and at least one positive solution to the problem (\ref{1}).
\end{teor}

Note that the system (\ref{1}) has a variational structure. In fact, it can be seen as a Hamiltonian
system, since, if we consider
\[
H(u, v,x) = \frac{v^{p+1}}{(p+1)\vert x\vert^\alpha}+ \frac{u^{q+1}}{(q+1)\vert x\vert^\beta} 
\]
then we have
\[
H_v( u, v,x) = \frac{v^p}{\vert x\vert^\alpha}\text{ and }H_u( u, v,x) = \frac{u^q}{\vert x\vert^\beta}.
\]

\begin{rem}
For such solutions, we prove an $L^\infty$ estimate of Brezis-Kato type and derive the regularity property of the solutions based on the results obtained in \cite{choi, ROS}.
\end{rem}

\begin{rem}
Related systems have been investigated by using other methods. We refer to the works \cite{serena2, serena1, EM} for systems involving different operators $(-\Delta)^{s}$ and $(-\Delta)^{t}$ in each one of equations. More generally, fractional systems have been studied with extension methods in \cite{serena3, serena4}.
\end{rem}

The rest of paper is organized into five sections. In Section 2 we briefly recall some definitions and facts dealing with fractional Sobolev spaces and comment some relationships and differences between operators $(-\Delta)^{s}$ and $\mathcal{A}^s$. In Section 3 we establish the functional setting in which the problem will be posed. In Section 4 we prove our main result, Theorem \ref{teo1}. Finally in Section 5 we shall establish the Brezis-Kato type result and study the regularity of solutions to (\ref{1}).

\section{Preliminaries}

In this section we briefly recall some definitions and facts related to two types of operators $\mathcal{A}^s$ and $(-\Delta)^s$.

$\bullet$ The spectral fractional Laplace operator: For $\Omega$ be a smooth bounded open subset of $\mathbb{R}^N$. The spectral fractional Laplace operator $\mathcal{A}^s$ is defined as follows. Let $\varphi_k$ be an eigenfunction of $-\Delta$ given by
\begin{equation}
\left\{\begin{array}{ccll}
-\Delta\varphi_k&=&\lambda_k\varphi_k & {\rm in} \ \ \Omega\; \\
\varphi_k&=& 0 & {\rm on} \ \ \partial\Omega
\end{array}\right.  ,
\end{equation}
where $\lambda_k$ is the corresponding eigenvalue of $\varphi_k,0<\lambda_1<\lambda_2\leq\lambda_3\leq\cdots\leq\lambda_k\rightarrow +\infty.$ Then, $\{\varphi_k\}_{k=1}^{\infty}$ is an orthonormal basic of $L^2(\Omega)$ satisfying 
\[
\int\limits_{\Omega}\varphi_j\varphi_kdx=\delta_{j,k}.
\]
We define the operator $\mathcal{A}^s$ for any $u\in C^\infty_0(\Omega)$ by

\begin{equation}\label{lapla}
\mathcal{A}^su=\sum_{k=1}^{\infty}\lambda_k^s\xi_k\varphi_k,
\end{equation}
where
\[
u=\sum_{k=1}^{\infty}\xi_k\varphi_k\text{  and  }\xi_k=\int\limits_{\Omega}u\varphi_kdx.
\]

$\bullet$ The restricted fractional Laplace operator: In this case we materialize the zero Dirichlet condition by restricting the operator to act only on functions that are zero outside $\Omega$. We will call the operator defined in such a way the restricted fractional Laplace operator. So defined, $(-\Delta)^s$ is a self-adjoint operator on $L^{2}(\Omega)$, with a discrete spectrum: we will denote by $\mu_{k}>0$, $k= 1, 2,\ldots$ its eigenvalues written in increasing order and repeated according to their multiplicity and we will denote by $\{\psi_{k}\}_k$ the corresponding set of eigenfunctions, normalized in $L^2(\Omega)$, where $\psi_{k}\in H^{2s}_0(\Omega)$. Eigenvalues $\mu_{k}$ (including multiplicities) satisfy
\[
0<\mu_{1}<\mu_{2}\leq\mu_{3}\leq\cdots\leq\mu_{k}\rightarrow +\infty.
\]

The spectral fractional Laplace operator $\mathcal{A}^s$ is related to (but different from) the restricted fractional Laplace operator $(-\Delta)^s$.

\begin{teor} The operators $(-\Delta)^s$ and $\mathcal{A}^s$ are not the same, since they have different eigenvalues and eigenfunctions. More precisely:
\begin{itemize}
    \item[(1)] the first eigenvalues of $(-\Delta)^s$ is strictly less than the one of $\mathcal{A}^s$.
    \item[(2)] the eigenfunctions of $(-\Delta)^s$ are only H\"{o}lder continuous up to the boundary, differently from the ones of $\mathcal{A}^s$ that are as smooth up the boundary as the boundary allows.
\end{itemize}
\end{teor}
\n {\bf Proof.} See \cite{rafaella}.\; \fim\\

The next theorem gives a relation between the spectral fractional Laplace operator $\mathcal{A}^s$ and the restricted fractional Laplace operator $(-\Delta)^s$.

\begin{teor}
For $u\in H^{s}(\mathbb{R}^N)$, $u\geq 0$ and $\hbox{supp}(u)\subset\overline{\Omega}$, the following relation holds in the sense of distributions:
\[
\mathcal{A}^su\geq (-\Delta)^su.
\]
If $u\neq 0$ then this inequality holds with strict sign.
\end{teor}
\n {\bf Proof.} See \cite{compare}.\; \fim\\

$\bullet$ Common notation. In the sequel we use $\mathcal{L}$ to refer to any of the two types of operators $\mathcal{A}^s$ or $(-\Delta)^s$, $0 < s < 1$. Each one is defined on a Hilbert space
\[
\Theta^{s}(\Omega)=\{u=\sum_{k=1}^{\infty}u_k\psi_{k}\in L^2(\Omega)\mid\sum_{k=1}^{\infty}\mu_{k}\vert u_k\vert^2<+\infty\}
\]
with values in its dual $\Theta^{s}(\Omega)'$. The Spectral Theorem allows to write $\mathcal{L}$ as
\[
\mathcal{L}u=\sum_{k=1}^{\infty}\mu_{k}u_k\psi_{k}
\]
for any $u\in \Theta^{s}(\Omega)$. Thus the inner product of $\Theta^{s}(\Omega)$ is given by
\[
\langle u,v\rangle_{\Theta^{s}(\Omega)}=\int\limits_{\Omega}\mathcal{L}^{1/2}u\mathcal{L}^{1/2}vdx=\int\limits_{\Omega}u\mathcal{L}vdx=\int\limits_{\Omega}v\mathcal{L}udx.
\]
We denote by $\Vert\cdot\Vert_{\Theta^{s}(\Omega)}$ the norm derived from this inner product. The notation in the formula copies the one just used for the second operator. When applied to the first one we put here $\psi_{k} = \varphi_k$, and $\mu_{k} = \lambda_{k}^s$. Note that $\Theta^{s}(\Omega)$ depends in principle on the type of operator and on the exponent $s$. It turns out that $\Theta^{s}(\Omega)$ independent of operator for each $s$, see \cite{sire}. We remark that $\Theta^{s}(\Omega)'$ can be described as the completion of the finite sums of the form
\[
f=\sum_{k=1}^{\infty}c_k\psi_{k} 
\]
with respect to the dual norm
\[
\Vert f\Vert_{\Theta^{s}(\Omega)'}=\sum_{k=1}^{\infty}\mu_{k}^{-1}\vert c_k\vert^2=\Vert\mathcal{L}^{-1/2}f\Vert_{L^2(\Omega)}^2=\int\limits_{\Omega}f\mathcal{L}^{-1}fdx
\]
and it is a space of distributions. Moreover, the operator $\mathcal{L}$ is an isomorphism between $\Theta^{s}(\Omega)$ and $\Theta^{s}(\Omega)'\simeq\Theta^{s}(\Omega)$, given by its action on the eigenfunctions. If $u,v\in\Theta^{s}(\Omega)$ and $f = \mathcal{L}u$ we have, after this isomorphism,
\[
\langle f, v\rangle_{\Theta^{s}(\Omega)'\times\Theta^{s}(\Omega)} = \langle u, v\rangle_{\Theta^{s}(\Omega)\times\Theta^{s}(\Omega)} =\sum_{k=1}^{\infty}\mu_{k}u_kv_k.
\]
If it also happens that $f\in L^{2}(\Omega)$, then clearly we get 
\[
\langle f, v\rangle_{\Theta^{s}(\Omega)'\times\Theta^{s}(\Omega)}=\int\limits_{\Omega}fv dx.
\]
We have $\mathcal{L}^{-1}:\Theta^{s}(\Omega)'\rightarrow\Theta^{s}(\Omega)$ can be written as
\[
\mathcal{L}^{-1}f(x)=\int\limits_{\Omega}G_{\Omega}(x,y)f(y)dy,
\]
where $G_{\Omega}$ is the Green function of operator $\mathcal{L}$ (see \cite{Green, ros176}). It is known that
\[
\Theta^s(\Omega)=\left\{
\begin{array}{llll}
L^2(\Omega) & {\rm if} \ \ s=0\\
H^{s}(\Omega)=H^{s}_0(\Omega) & {\rm if} \ \ s\in(0,\frac{1}{2})\\
H^{\frac{1}{2}}_{00}(\Omega) & {\rm if} \ \ s=\frac{1}{2}\\
H^{s}_0(\Omega) & {\rm if} \ \ s\in(\frac{1}{2},1]\\
H^{s}(\Omega)\cap H^{1}_0(\Omega) & {\rm if} \ \ s\in(1,2]
\end{array},
\right.
\]
where $H^{\frac{1}{2}}_{00}(\Omega):=\{u\in H^{1/2}(\Omega)\mid\int_{\Omega}\frac{u^2(x)}{d(x)}dx<+\infty\}.$\\

Observe that the injection $\Theta^s(\Omega)\hookrightarrow H^{s}(\Omega)$ is continuous. By the Sobolev imbedding theorem we therefore have continuous imbeddings 
\begin{equation}\label{1.1}
\Theta^s(\Omega)\subset L^{p+1}(\Omega)
\end{equation} 
if $p+1\leq\frac{2N}{N-2s}$ and these imbedding are compact if $p+1<\frac{2N}{N-2s}$ for $0<s<2N$. Also, we have compact imbedding  $\Theta^s(\Omega)\subset C(\Omega)$, if 
\[
\frac{s}{N}>\frac{1}{2}.
\]
For $0<r<2$ we have $\mathcal{L}:\Theta^{r}(\Omega)\rightarrow\Theta^{r-2s}(\Omega)$ is an isomorphism (see \cite{vander}). Now, we use Holder's inequality to obtain
\[
\int\limits_{\Omega}\frac{u^{q+1}}{\vert x\vert^{\beta}}dx\leq\left(\int\limits_{\Omega}u^{\gamma}dx\right)^{(q+1)/\gamma}\left(\int\limits_{\Omega}\vert x\vert^{-\beta\gamma/(\gamma-(q+1))}dx\right)^{(\gamma-(q+1))/\gamma}\leq C\left(\int\limits_{\Omega}u^{\gamma}dx\right)^{(q+1)/\gamma}
\]
if 
\[
\frac{\beta\gamma}{\gamma-(q+1)}<N.
\]
That is
\[
N(q+1)<(N-\beta)\gamma.
\]

\begin{propo}\label{P2.1} In the hypothesis of Theorem \ref{teo1}, we have 
\[\Theta^{t}(\Omega)\times\Theta^{2s-t}(\Omega)\hookrightarrow L^{q+1}(\Omega,\vert x\vert^{-\beta})\times L^{p+1}(\Omega,\vert x\vert^{-\alpha})
\]
is compact.
\end{propo}

\section{The variational formulation of system (\ref{1})}
We will the proof for the spectral fractional Laplace operator $\mathcal{A}^s$. Similarly, follows the results for the restricted fractional Laplace operator $(-\Delta)^s$, changing the corresponding space. 

The existence result follows by applying the proof of \cite{rossi} for the case $s = 1$ with only minor modifications.

We define the product Hilbert spaces
\[
E^{t}(\Omega)=\Theta^{t}(\Omega)\times\Theta^{2s-t}(\Omega),\text{   }0<t<2s
\]
where your inner product is given by
\[
\langle(u_1,v_1),(u_2,v_2)\rangle_{E^{t}(\Omega)}=\langle\mathcal{A}^{t/2}u_1,\mathcal{A}^{t/2}u_2\rangle_{L^{2}(\Omega)}+\langle\mathcal{A}^{s-t/2}v_1,\mathcal{A}^{s-t/2}v_2\rangle_{L^{2}(\Omega)}.
\]

We denote by $\Vert\cdot\Vert_{E}$ the norm derived from this inner product, i.e,

\[
\Vert(u,v)\Vert_{E}=\left(\Vert u\Vert^2_{\Theta^{t}}+\Vert v\Vert^2_{\Theta^{2s-t}}\right)^{\frac{1}{2}}.
\]

We also have $\mathcal{A}^s:\Theta^t(\Omega)\rightarrow\Theta^{t-2s}(\Omega)$ is an isomorphism, see \cite{vander}. Hence
$$\left(\begin{array}{ccllrr}
0& \mathcal{A}^{s}  \\
\mathcal{A}^{s} & 0 
\end{array}\right): E^t(\Omega)\rightarrow\Theta^{-t}\times\Theta^{t-2s}(\Omega)=E^t(\Omega)'$$
is an isometry. We consider the Lagrangian
\begin{eqnarray*}
\mathcal{J}(u,v)&=&\int\limits_{\Omega}\mathcal{A}^{s/2}u\mathcal{A}^{s/2}vdx-\int\limits_{\Omega}H(u(x),v(x),x)dx,
\end{eqnarray*}
i.e., a strongly indefinite functional. The o Hamiltonian is given by
\begin{equation}\label{2.2.}
\mathcal{H}(u,v)=\int\limits_{\Omega}H(u(x),v(x),x))dx.
\end{equation}
The quadratic part can again be written as
\[
Q(u,v)=\frac{1}{2}\langle L(u,v),(u,v)\rangle_{E^{\alpha}(\Omega)}=\int\limits_{\Omega}\mathcal{A}^{t/2}u\mathcal{A}^{s-t/2}vdx=\int\limits_{\Omega}\mathcal{A}^{s/2}u\mathcal{A}^{s/2}vdx,
\]
where
\begin{equation}\label{1.8}
L=\left(\begin{array}{ccllrr}
0& \mathcal{A}^{s-t}  \\
\mathcal{A}^{t-s} & 0 
\end{array}\right)
\end{equation}
is bounded and self-adjoint. Introducing the "diagonals"
\[
E^+=\{(u,\mathcal{A}^{t-s}u):u\in\Theta^t(\Omega)\}\text{ and }E^-=\{(u,-\mathcal{A}^{t-s}u):u\in\Theta^t(\Omega)\}
\]
we have
\[
E^{t}(\Omega)=E^+\oplus E^-.
\]

\begin{propo} The functional $\mathcal{H}$ defined in (\ref{2.2.}) is of class $C^1$ and its derivative
is given by
\begin{equation}\label{1.29}
\mathcal{H}'(u,v)(\varphi,\psi)=\int\limits_\Omega\frac{\partial H}{\partial u}(u,v,x)\varphi+ \frac{\partial H}{\partial v}(u,v,x)\psi dx
\end{equation}
for all $(u, v),(\varphi,\psi)\in E^t(\Omega)$. Moreover $\mathcal{H}' : E^t(\Omega)\rightarrow E^t(\Omega)$ is a compact operator.
\end{propo}

{\bf Proof.} The expression given at the right-hand side of (\ref{1.29}) is well defined. In fact, we have
\[
\int\limits_\Omega\left\vert\frac{\partial H}{\partial u}(u,v,x)\varphi\right\vert dx=\int\limits_\Omega\left(\frac{\vert u\vert^q}{\vert x\vert^{\beta}}\right)\vert\varphi\vert dx.
\]
Using Holder's inequality and Proposition \ref{P2.1} we have
\[
\int\limits_\Omega\left\vert\frac{\partial H}{\partial u}(u,v,x)\varphi\right\vert dx\leq C\left(\Vert u\Vert_{\Theta^t}^q\right)\Vert\varphi\Vert_{\Theta^t}.
\]
In a similar way we obtain an inequality for the derivative with respect to $v$. Thus $\mathcal{H}'(u, v)$ is well defined and bounded in $E^t(\Omega)$.

Next, usual arguments give that $\mathcal{H}$ is Fr\'{e}chet differentiable, $\mathcal{H}'$ is continuous
and, as a consequence of the Proposition \ref{P2.1}, $\mathcal{H}'$ is also compact. See \cite{Rabinowitz} for example. \; \fim\\

\begin{propo}\label{imer} In the hypothesis of Theorem \ref{teo1}, if $(u,v)\in E^t(\Omega)$, we have $\frac{u^{q-1}}{\vert x\vert^{\beta}}\in L^d(\Omega)$ and $\frac{v^{p-1}}{\vert x\vert^{\alpha}}\in L^c(\Omega)$ for every
\[
1<c<\frac{2N}{(p-1)(N-(4s-2t))+2\alpha}\text{ and }1<d<\frac{2N}{(q-1)(N-2t)+2\beta}
\]
whenever $p,q\geq 1$ and we have $\frac{u^{q/2}}{\vert x\vert^{\beta}}\in L^d(\Omega)$ and $\frac{v^{p/2}}{\vert x\vert^{\alpha}}\in L^c(\Omega)$ for every
\[
1<c<\frac{4N}{p(N-(4s-2t))+4\alpha}\text{ and }1<d<\frac{4N}{q(N-2t)+4\beta}
\]
whenever $p,q\in(0,1)$.
\end{propo}

\section{Proof of Theorem \ref{teo1}}

First, we prove that there exist infinitely many solutions to (\ref{1}). To this end, we present an abstract theorem from critical point theory from \cite{23} (see also \cite{7}) that provides us with infinitely many critical points of a functional. Next, we prove that this abstract result can be applied to our functional setting stated in the previous section.

Let $E$ be a Hilbert space with inner product $\langle\cdot,\cdot\rangle_E$. Assume that $E$ has a splitting $E= X\oplus Y$ where $X$ and $Y$ are both infinite dimensional subspaces. Assume there exists a sequence of finite dimensional subspaces $X_n\subset X$, $Y_n\subset Y$, $E_n = X_n \oplus Y_n$ such that $\overline{\cup_{n=1}^{\infty}E_n} = E$. Let $T : E \rightarrow E$ be a linear bounded invertible operator.

Let $\mathcal{J}\in C^1(E,\mathbb{R})$. Instead of the usual Palais-Smale condition we will require that the functional $\mathcal{J}$ satisfies the so-called $(PS)^\ast$ conditions with respect to $E_n$, i.e., any sequence $z_k\in E_{n_k}$ with $n_k\rightarrow\infty$ as $k\rightarrow\infty$, satisfying $\mathcal{J}\vert_{E_{n_k}}'(z_k)\rightarrow 0$ and $\mathcal{J}(z_k)\rightarrow c$ has a subsequence that converges in $E$.

Then we define the basic sets over which the linking process will take place. For $\rho > 0$ we define
\[
S = S_\rho = \{y\in Y: \Vert y\Vert_E = \rho\}
\]
and for some fixed $y_1\in Y$ with $\Vert y_1\Vert_E = 1$ and subspaces $\mathcal{X}_1$ and $\mathcal{X}_2$, we consider 
\[
X \oplus [y_1] = \mathcal{X}_1\oplus\mathcal{X}_2.
\]
Without loss of generality we may assume that $y_1\in \mathcal{X}_2$. Next, we define for $M, \sigma > 0$
\[
D = D_{M,\rho} = \{x_1 + x_2\in\mathcal{X}_1\oplus\mathcal{X}_2: \Vert x_1\Vert_E \leq M, \Vert x_2\Vert_E\leq\sigma\}.
\]
Now we can state our abstract critical point result whose proof can be found in \cite{23}.

\begin{teor}\label{teo3.1} Let $\mathcal{J}\in C^1(E,\R)$ be an even functional satisfying the $(PS)^\ast$ condition with respect to $E_n$. Assume that $T : E_n\rightarrow E_n$, for $n$ large. Let $\rho > 0$ and $\sigma > 0$ be such that $\sigma\Vert T y_1\Vert_E > \rho$. Assume that there are constants $\alpha\leq\beta$ such that
\[
\inf_{S\cap E_n}\mathcal{J}\geq\alpha,\text{  }\sup_{T(\partial D\cap E_n)}\mathcal{J} <\alpha\text{ and }\sup_{T(D\cap E_n)}\mathcal{J}\leq\beta
\]
for all $n$ large. Then $\mathcal{J}$ has a critical value $c\in [\alpha,\beta]$.
\end{teor}

Next, we show how the functional setting introduced in Section 3 can be used to apply Theorem \ref{teo3.1}. Let
\[
E_n = [\varphi_1,\ldots,\varphi_n] \times [\varphi_1,\ldots,\varphi_n].
\]
It is easy to see that $\overline{\cup_{n=1}^{\infty}E_n} = E^t(\Omega)$. Next, we prove that $\mathcal{J}$ satisfies the $(PS)^\ast$ condition with respect to the family $E_n$.

\begin{lema}\label{lema3.1} The functional $\mathcal{J}$ satisfies the $(PS)^\ast$ condition with respect to $E_n$. 
\end{lema}
{\bf Proof.} Let $(z_k)_{k\geq 1}=(u_k,v_k)_{k\geq 1}\subset E_{n_k}$ be a sequence such that
\begin{equation}\label{3.1}
\mathcal{J}(z_k)\rightarrow c,\text{ and }\mathcal{J}\vert_{E_{n_k}}'(z_k)\rightarrow 0.
\end{equation}
Let us first prove that (\ref{3.1}) implies that $(z_k)$ is bounded in $E^t(\Omega)$. From (\ref{3.1}) it follows that there exists a sequence $(\varepsilon_k)$ converging to $0$ such that
\begin{equation}\label{3.2}
\vert\mathcal{J}'(z_k)w\vert\leq\varepsilon_k\Vert w\Vert_E,\forall w\in E_{n_k}.
\end{equation}
Let us take
\[
w_k=((w_k)_1,(w_k)_2)\frac{(q+1)(p+1)}{p+q+2}\left(\frac{1}{q+1}u_k,\frac{1}{p+1}v_k\right),\text{ where }z_k=(u_k,v_k).
\]
Now, using (\ref{3.1}) and (\ref{3.2}), for $k$ large
\begin{eqnarray*}
c + 1+ \varepsilon_k \Vert w_k \Vert_{E} &\geq &  \mathcal{J}(z_k) - \mathcal{J}'(z_k)w_k \\
&= &\int\limits_{\Omega}\mathcal{A}^{t/2}u_k\mathcal{A}^{s-t/2}v_kdx- \int\limits_{\Omega} H(u_k,v_k,x)dx-\int\limits_{\Omega}\mathcal{A}^{t/2}u_k\mathcal{A}^{s-t/2}(w_k)_2dx\\
&- &\int\limits_{\Omega}\mathcal{A}^{t/2}(w_k)_1\mathcal{A}^{s-t/2}v_kdx+ \int\limits_{\Omega}H_u(u_k,v_k,x)(w_k)_1 +\int\limits_{\Omega} H_v(u_k,v_k,x)(w_k)_2 dx\\
&= & -\frac{(1-pq)}{p+q+2} \int\limits_{\Omega}H(u_k,v_k,x)dx.
\end{eqnarray*}
Now, by (\ref{1.4}) we get $pq>1$ and hence we obtain 
\[
C(1 + \Vert z_k\Vert_E)\geq \int\limits_{\Omega}H(u_k,v_k,x)dx.
\]
Therefore,
\begin{equation}\label{3.3}
\int\limits_{\Omega}\frac{\vert u_k\vert^{q+1}}{\vert x\vert^\beta} + \frac{\vert v_k\vert^{p+1}}{\vert x\vert^{\alpha}}dx \leq C(1 + \Vert u_k\Vert_{\Theta^t}+ \Vert v_k\Vert_{\Theta^{2s-t}}).
\end{equation}
Next, let us consider $w= (\phi,0)$ with $\phi\in\Theta^t_{n_k}(\Omega)$. Then from (\ref{3.2})
\begin{equation}\label{2.6}
\left\vert\int\limits_{\Omega}\mathcal{A}^{t/2}\phi\mathcal{A}^{s-t/2}v_kdx\right\vert\leq\int\limits_{\Omega}\frac{\vert u_k\vert^q}{\vert x\vert^\beta}\vert\phi\vert dx+\varepsilon_k\Vert\phi\Vert_{\Theta^t}.
\end{equation}
Now, using Holder's inequality,
\[
\int\limits_{\Omega}\frac{\vert u_k\vert^{q}}{\vert x\vert^\beta}\vert\phi\vert dx\leq \Vert u_k\Vert^{q}_{L^{q+1}(\Omega,\vert x\vert^{-\beta})}\Vert\phi\Vert_{L^{q+1}(\Omega,\vert x\vert^{-\beta})}.
\]
So, by Proposition \ref{P2.1}, we get
\[
\vert\langle\mathcal{A}^{t/2}\phi,\mathcal{A}^{s-t/2}v_k\rangle\vert\leq C\Vert\phi\Vert_{\Theta^t}\left(\Vert u_k\Vert^q_{L^{q+1}(\Omega,\vert x\vert^{-\beta})} + 1\right).
\]
By duality ($\mathcal{A}^{t/2}$ is an isometry between $\Theta^t(\Omega)$ and $L^2(\Omega)$) we get
\begin{equation}\label{3.4}
\Vert v_k\Vert_{\Theta^{2s-t}}\leq C\left(\Vert u_k\Vert^q_{L^{q+1}(\Omega,\vert x\vert^{-\beta})} + 1\right).
\end{equation}
By an analogous reasoning
\begin{equation}\label{3.5}
\Vert u_k\Vert_{\Theta^{t}}\leq C\left(\Vert v_k\Vert^p_{L^{p+1}(\Omega,\vert x\vert^{-\alpha})} + 1\right).
\end{equation}
Now combining (\ref{3.3}), (\ref{3.4}) and (\ref{3.5}), we obtain
\[
\Vert u_k\Vert_{\Theta^{t}}+\Vert v_k\Vert_{\Theta^{2s-t}}\leq C\left(\Vert u_k\Vert^{q/(q+1)}_{\Theta^{t}} + \Vert v_k\Vert^{p/(p+1)}_{\Theta^{2s-t}}+1\right).
\]
Since all the involved exponents are less than one, we conclude that $z_k$ in bounded. Now, by compactness and the invertibility of $L$ we can extract a subsequence of $z_k$ that converges in $E^t(\Omega)$. Indeed, we can take a subsequence $z_{k_j}$ that converges weakly in $E^t(\Omega)$, as $\mathcal{H}$ is compact, it follows that $\mathcal{H}'(z_{k_j})$ converges strongly in $E^t(\Omega)$. Hence, using the fact that $\mathcal{J}(z_{k_j})\rightarrow 0$ strongly and the invertibility of $L$, the result follows.\; \fim\\

Now we define the splitting of $E_n$. Fix $k\in\mathbb{N}$ and for $n\geq k$ let
\begin{equation}
X_n = (E^-_1\oplus\ldots\oplus E^-_n) \oplus (E^+_1\oplus\ldots E^+_{k-1})\text{ and }Y_n = (E^+_k\oplus\ldots\oplus E^+_n),
\end{equation}
where $E^+_j= [(\varphi_j,\mathcal{A}^{t-s}\varphi_j)]$ and $E^-_j= [(\varphi_j,-\mathcal{A}^{t-s}\varphi_j)]$. We have $E_n =X_\oplus Y_n$.

\begin{lema}\label{lema3.2} There exist $\alpha_k > 0$ and $\rho_k > 0$ independent of $n$ such that for all $\geq k$
\[
\inf_{z\in S_{\rho_k}\cap Y_n}\mathcal{J}(z)\geq\alpha_k
\]
where $S_{\rho_k}= \{y\in E^+: \Vert y\Vert =\rho_k\}$. Moreover, $\alpha_k\rightarrow\infty$ as $k\rightarrow\infty$.
\end{lema}
{\bf Proof.} We first recall that by Proposition \ref{P2.1}, $\Theta^t(\Omega)$ is embedded in $L^\gamma(\Omega,\vert x\vert^{-\varrho})$ for any $\gamma$ such that
\[
\gamma\leq\frac{2(N-\varrho)}{N-2t}.
\]
Hence, there exists $a = a(\gamma)$ such that
\[
\Vert u\Vert_{L^\gamma(\Omega,\vert x \vert^{-\varrho})}\leq a\Vert u\Vert_{\Theta^t}\text{ for all }u\in\Theta^t(\Omega).
\]
Also for $z\in E^+_k\oplus\ldots\oplus E^+_j\oplus\ldots$ we have
\[
\Vert z\Vert_E\geq\lambda_k^{s\min\{t,2s-t\}}\Vert z\Vert_{L^2}
\]
with $\lambda_k\rightarrow\infty$ as $k\rightarrow\infty$.

Now consider $z = (u, v)\in Y_n$. For a constant $a$ independent of $n$, we observe that there exists $\kappa > 0$ such that
\[
\Vert u\Vert^{q+1}_{L^{q+1}(\Omega,|x|^{-\beta})}\leq \Vert u\Vert^{2/\kappa}_{L^{2}}\Vert u\Vert^{q+1-2/\kappa}_{L^{\gamma}(\Omega,|x|^{-\varrho})}\leq \frac{a}{\lambda_k^{s\min\{t,2s-t\}(2/\kappa)}}\Vert u\Vert^{q+1}_E.
\]
Analogously, we obtain
\[
\Vert v\Vert^{p+1}_{L^{p+1}(\Omega,|x|^{-\alpha})}\leq \frac{a}{\lambda_k^{s\min\{t,2s-t\}(2/\theta)}}\Vert v\Vert^{p+1}_E
\]
for some $\theta> 0$. Then for $z = (u, v)$ we have
\[
\mathcal{J}(z)\geq \Vert z\Vert^2_E- C\left(\frac{a}{\lambda_k^{s\min\{t,2s-t\}\min\{2/\kappa,2/\theta\}}}\max \left\{\Vert z\Vert^{q+1}_E, \Vert z\Vert^{p+1}_E\right\}+1\right).
\]
Then we choose 
\[
\rho^{\max\{p+1,q+1\}}_k=\lambda_k^{s\min\{t,2s-t\}\min\{2/\kappa,2/\theta\}}
\]
and observe that $\rho_k\rightarrow\infty$ as $k\rightarrow\infty$.
Therefore, for $z\in S_{\rho_k}\cap Y_n$ we find
\begin{equation}\label{3.7}
\mathcal{J}(z)\geq \rho_k^2- C. 
\end{equation}
Defining $\alpha_k$ as the right hand side of (\ref{3.7}) and noting that both $\rho_k$ and $\alpha_k$ are independent of $n\geq k$ we complete the proof of the Lemma.\; \fim\\

Next we define, for $z = (u, v)\in E^t(\Omega)$
\[
T_\sigma(z) = (\sigma^{\mu-1}u,\sigma^{\nu-1}v) 
\]
where $\mu$ and $\nu$ are such that
\[
\mu+\nu < \min\{\mu(p + 1), \nu(q + 1)\}.
\]
In particular, since $pq > 1$ by (\ref{1.4}), we take $\mu = q + 1$ and $\nu = q + 1$ for which,
\[
(p + 1) + (q + 1) < (p + 1)(q + 1).
\]

\begin{lema}\label{lema3.3} There exist $B_k > 0$, $\sigma_k$ and $M_k > 0$ independent of $n$ such that for all $n \geq k$ they satisfy $\sigma_k >\rho_k$,
\[
\sup_{T_{\sigma_k}(\partial D\cap E_n)}\mathcal{J}\leq 0\text{ and }\sup_{T_{\sigma_k}(D\cap E_n)}\mathcal{J}\leq B_k
\]
where
\[
D = \{z\in E^-\oplus E^+_1\oplus\ldots\oplus E^+_k:\Vert z^-\Vert\leq M_k, \Vert z^+\Vert\leq\sigma_k\}.
\]
\end{lema}
{\bf Proof.} Let us consider $z = T_\sigma(u, v)$ with $(u, v)\in D$. Then we can write
\[
z = (\sigma^{\mu-1}u^+,\sigma^{\nu-1}v^+)+(\sigma^{\mu-1}u^-,\sigma^{\nu-1}v^-).
\]
Using the definition of the spaces $E^+$ and $E^-$ we have
\[
\int\limits_{\Omega}\mathcal{A}^{t/2}u\mathcal{A}^{s-t/2}vdx = \sigma^{\mu+\nu-2}(\Vert z^+\Vert^2 - \Vert z^-\Vert^2).
\]
On the other hand, we have
\[
\int\limits_{\Omega}H(z,x) dx=\int\limits_{\Omega}\left(\sigma^{(q+1)(\mu-1)}\frac{\vert u^++u^-\vert^{q+1}}{\vert x\vert^\beta}+\sigma^{(p+1)(\nu-1)}\frac{\vert v^++v^-\vert^{p+1}}{\vert x\vert^\alpha}\right)dx.
\]
The functions $u^+$ and $u^-$ can be written as
\[
u^+ =\sum_{i=1}^{k}\theta_i\varphi_i\text{ and }u^- =\sum_{i=1}^{k}\gamma_i\varphi_i+\tilde{u}^-,
\]
where $\tilde{u}^-$ is orthogonal to $\varphi_i$, $i = 1,\ldots, k$ in $L^2(\Omega)$. Using Holder's inequality we get
\begin{eqnarray*}
\sum_{i=1}^{k}\lambda^{t-2}(\theta_i^2+\theta_i\gamma_i)&=&\langle u^++u^-,\mathcal{A}^{t-s}u^+\rangle\leq\Vert u^+ + u^-\Vert_{L^{q+1}(\Omega,|x|^{-\beta})}\Vert\mathcal{A}^{t-s}u^+\Vert_{L^{\frac{q+1}{q}}(\Omega,|x|^{\frac{\beta}{q}})}\\
&\leq &C(\Omega,\beta)\Vert u^+ + u^-\Vert_{L^{q+1}(\Omega,|x|^{-\beta})}\Vert\mathcal{A}^{t-s}u^+\Vert_{L^2}.
\end{eqnarray*}
Then by the definition of $\mathcal{A}^{t-s}$, there exists a constant $C_k$ such that
\begin{equation}\label{3.9}
\sum_{i=1}^{k}\lambda_i^{t-s}(\theta_i^2+\theta_i\gamma_i)\leq  C_k\Vert u^+ + u^-\Vert_{L^{q+1}(\Omega,|x|^{-\beta})}\Vert u^+\Vert_{L^2}.
\end{equation}
In a similar way, using that $v^+=\mathcal{A}^{t-s}u^+$ and $v^- = -\mathcal{A}^{t-s}u^-$ we have there exists a constant $C_k$ such that
\begin{equation}\label{3.10}
\sum_{i=1}^{k}\lambda_i^{t-s}(\theta_i^2+\theta_i\gamma_i)\leq  C_k\Vert v^+ + v^-\Vert_{L^{p+1}(\Omega,|x|^{-\alpha})}\Vert v^+\Vert_{L^2}.
\end{equation}
Depending on the sign $\sum_{i=1}^{k}\lambda_i^{t-s}\theta_i\gamma_i$ we use (\ref{3.9}) or (\ref{3.10}) to conclude that
\[
\Vert u^+\Vert_{L^2}\leq C_k\Vert u^+ + u^-\Vert_{L^{q+1}(\Omega,|x|^{-\beta})}
\]
or
\[
\Vert u^+\Vert_{L^2}\leq C_k\Vert v^+ + v^-\Vert_{L^{p+1}(\Omega,|x|^{-\alpha})}.
\]
Hence,
\[
\mathcal{J}(z)\leq\sigma^{\mu+\nu-2}(\Vert z^+\Vert^2 - \Vert z^-\Vert^2) -C_k\sigma^{(q+1)(\mu-1)}\Vert u^+\Vert^{q+1}_{L^2}
\]
or
\[
\mathcal{J}(z)\leq\sigma^{\mu+\nu-2}(\Vert z^+\Vert^2 - \Vert z^-\Vert^2) -C_k\sigma^{(p+1)(\nu-1)}\Vert u^+\Vert^{p+1}_{L^2}.
\]
Thus we may choose $\Vert z^+\Vert_E=\sigma_k$ large enough in order to obtain $\sigma_k >\rho_k$ and, by the condition on $\mu$ and $\nu$, $\mathcal{J} (z)\leq 0$.
Taking $\Vert z^+\Vert\leq\sigma_k$ and $\Vert z^-\Vert = M_k$, we get
\[
\mathcal{J}(z)\leq\sigma_k^{\mu+\nu-2}(\sigma_k^2-M_k^2) 
\]
and then choosing $M_k$ large enough we find that
\[
\mathcal{J}(z)\leq 0.
\]
In this way we have finished with the proof of the first part of Lemma \ref{lema3.3}. 

Next we choose $B_k$ large so that the second inequality holds.\; \fim\\

{\bf Proof of Theorem \ref{teo1}.} Existence of infinitely many solutions.

For $k\geq 1$, Lemmas \ref{lema3.2} and \ref{lema3.3} allow us to use Theorem \ref{teo3.1}. As a consequence the functional $\mathcal{J}$ has a critical value $c_k\in [\alpha_k, B_k]$. Since $\alpha_k\rightarrow\infty$ we get infinitely many critical values of $\mathcal{J}$. Therefore we have infinitely many solutions of (\ref{1}).\; \fim\\

Now we turn our attention to the existence of a positive solution to (\ref{1}). We use ideas from \cite{E4} under the functional setting of Section 3. We start by redefining the Hamiltonian. Let us define $\tilde{H}:\R\times\R\times\partial\Omega\rightarrow\R$ by
\[
\tilde{H}(u,v,x)=
\left\{
\begin{array}{llll}
H(u,v,x) & {\rm if} \ \ u,v\geq 0\\
H(0,v,x) & {\rm if} \ \ u\leq 0,v\geq 0\\
H(u,0,x) & {\rm if} \ \ u\geq 0,v\leq 0\\
0 & {\rm if} \ \ u,v\leq 0
\end{array}
\right..
\]
We observe that if $(u, v)$ is a nontrivial solution of
\begin{equation}\label{3.12}
\left\{
\begin{array}{llll}
\mathcal{A}^s u= \tilde{H}_v(u,v,x)& {\rm in} \ \ \Omega\\
\mathcal{A}^s v = \tilde{H}_u(u,v,x)& {\rm in} \ \ \Omega\\
u= v=0 & {\rm on} \ \ \partial\Omega
\end{array}
\right.,
\end{equation}
then by the maximum principle (see \cite{T}), we have $u$ and $v$ are strictly positive in $\Omega$. Hence $(u, v)$ is a positive solution of (\ref{1}).

To find a nontrivial solution of (\ref{3.12}) we want to apply the results of Section 4. By our assumptions, the new Hamiltonian $\tilde{H}$ is regular.

We have to adapt the proof of Theorem \ref{teo1} to the functional $\mathcal{J}$ with the Hamiltonian $H$ replaced by $\tilde{H}$. We observe that the proof of the Palais-Smale condition and the geometric conditions follows as before with some minor modifications, see \cite{E4} for the details.\\

{\bf Proof of Theorem \ref{teo1}.} Existence of a positive solution.
As a consequence of the previous results the modified functional $\mathcal{J}$ (with the modified
Hamiltonian $\tilde{H}$ instead of $H$) has a critical value $c\neq 0$. Hence, by the maximum principle,
we obtain a positive solution of (\ref{1}).\; \fim

\section{Regularity of solutions of systems (\ref{1})}

Next we prove the $L^\infty$ estimate of Brezis-Kato type.

\begin{propo}\label{prop}
Let $(u,v)$ be a solution of the problem (\ref{1}). In the hypothesis of Theorem \ref{teo1}, we have $(u, v) \in L^{\infty}(\Omega) \times L^{\infty}(\Omega)$ and, moreover, if $\alpha,\beta<2s$ we have $(u,v) \in C^{\sigma}(\R^N) \times C^{\sigma}(\R^N)$ for some $\sigma \in (0, 1)$.
\end{propo}
\n {\bf Proof.} In the hypothesis of Theorem \ref{teo1}, we take $0< p \leq 1$ and $q>1$. We rewrite the problem (\ref{1}) as follows
\begin{equation}
\left\{
\begin{array}{llll}
\mathcal{L} u=a(x)v^{\frac{p}{2}} & {\rm in} \ \ \Omega\\
\mathcal{L} v =b(x)u & {\rm in} \ \ \Omega\\
u= v=0 & {\rm on} \ \ \Sigma
\end{array}
\right.,
\end{equation}
where $a(x) = \dfrac{v^{\frac{p}{2}}}{\vert x\vert^\alpha}$ and $b(x) = \dfrac{u^{q-1}}{\vert x\vert^\beta}$. By Proposition \ref{imer}, we have $a \in L^{c}(\Omega)$ and $b \in L^{d}(\Omega)$ for every
\[
1<c<\frac{4N}{p(N-(4s-2t))+4\alpha}\text{ and }1<d<\frac{2N}{(q-1)(N-2t)+2\beta}.
\]
Thus, for each fixed $\varepsilon > 0$, we can construct functions $q_{\varepsilon} \in L^{c}(\Omega)$, $f_{\varepsilon} \in L^{\infty}(\Omega)$ and a constant $K_{\varepsilon} > 0$ such that

\[
a(x) v(x)^{\frac{p}{2}} = q_{\varepsilon}(x)v(x)^{\frac{p}{2}} + f_{\varepsilon}(x)
\]
and

\[
\Vert q_{\varepsilon} \Vert_{L^{c}} < \varepsilon,\ \ \Vert f_{\varepsilon} \Vert_{L^{\infty}} < K_{\varepsilon}\, .
\]
In fact, consider the set

\[
\Omega_k = \{x \in\Omega: \vert a(x)\vert < k\}\, ,
\]
where $k$ is chosen such that

\[
\int\limits_{\Omega_k^c} \vert a(x)\vert^{c}dx < \frac{1}{2}\varepsilon^{c}\, .
\]
This condition is clearly satisfied for $k = k_\varepsilon$ large enough.

We now write

\begin{equation}
q_{\varepsilon}(x)=\left\{
\begin{array}{llll}
\frac{1}{m}a(x) & {\rm for} \ \ x \in \Omega_{k_\varepsilon}\\
a(x) & {\rm for} \ \ x \in \Omega_{k_\varepsilon}^c\\
\end{array}
\right.
\end{equation}
and

\[
f_{\varepsilon}(x) = \left(a(x) - q_{\varepsilon}(x)\right)v(x)^{\frac{p}{2}}\, .
\]
Then,

\begin{eqnarray*}
\int\limits_{\Omega}\vert q_{\varepsilon}(x)\vert^{c}dx&=&\int\limits_{\Omega_{k_\varepsilon}}\vert q_{\varepsilon}(x)\vert^{c}dx+\int\limits_{\Omega_{k_\varepsilon}^c}\vert q_{\varepsilon}(x)\vert^{c}dx\\
&=&\left(\frac{1}{m}\right)^{c}\int\limits_{\Omega_{k_\varepsilon}}\vert a(x)\vert^{c}dx+\int\limits_{\Omega_{k_\varepsilon}^c}\vert a(x)\vert^{c}dx\\
&<&\left(\frac{1}{m}\right)^{c}\int\limits_{\Omega_{k_\varepsilon}}\vert a(x)\vert^{c}dx+\frac{1}{2}\varepsilon^{c}\, .
\end{eqnarray*}
So, for $m = m_\varepsilon > \left(\frac{2^{\frac{1}{c}}}{\varepsilon}\right) \Vert a\Vert_{L^{c}}$, we get

\[
\Vert q_{\varepsilon}\Vert_{L^{c}} < \varepsilon\, .
\]
Note also that $f_{\varepsilon}(x) = 0$ for all $x \in \Omega_{k_\varepsilon}^c$ and, for this choice of $m$,

\[
f_{\varepsilon}(x) = \left( 1-\frac{1}{m_\varepsilon}\right) a(x)^2 \leq \left( 1-\frac{1}{m_\varepsilon}\right) k_\varepsilon^2
\]
for all $x \in \Omega_{k_\varepsilon}$. Therefore,

\[
\Vert f_{\varepsilon} \Vert_{L^{\infty}} \leq \left( 1-\frac{1}{m_\varepsilon}\right) k_\varepsilon^2 := K_\varepsilon\, .
\]
On the other hand, we have

\[
v(x) = \mathcal{L}^{-1}(bu)(x)\, ,
\]
where $b\in L^{d}(\Omega)$. Hence,

\[
u(x) = \mathcal{L}^{-1}\left[q_{\varepsilon}(x)(\mathcal{L}^{-1}(bu)(x))^{\frac{p}{2}}\right] + \mathcal{L}^{-1}f_{\varepsilon}(x)\, .
\]

By Lemma 2.1 of \cite{choi} for $\mathcal{L}=\mathcal{A}^s$ and Proposition 1.4 of \cite{ROS} for $\mathcal{L}=(-\Delta)^s$, the claims $(ii)$ and $(iv)$ below follow readily and, by using H\"{o}lder's inequality, we also get the claims $(i)$ and $(iii)$. Precisely, for fixed $\gamma > 1$, we have:

 \begin{itemize}
    \item[(i)] The map $w \rightarrow b(x)w$ is bounded from $L^\gamma(\Omega)$ to $L^\tau(\Omega)$ for
\[
\frac{1}{\tau} = \frac{1}{d} + \frac{1}{\gamma};
\]
    \item[(ii)] For $\theta$ given by
\[
2s=N\left(\frac{1}{\tau}-\frac{2}{p \theta}\right),
\]
there exists a constant $C > 0$, depending on $\tau$ and $\theta$, such that

\[
\Vert (\mathcal{L}w)^{\frac{p}{2}}\Vert_{L^\theta}\leq C\Vert w\Vert_{L^\beta}^{\frac{p}{2}}
\]
for all $w \in L^\tau(\Omega)$;

    \item[(iii)] The map $ w \rightarrow q_\varepsilon(x)w$ is bounded from $L^\theta(\Omega)$ to $L^\eta(\Omega)$ with norm given by $\Vert q_\varepsilon \Vert_{L^{c}}$, where $\theta \geq 1$ and $\eta$ satisfies
\[
\frac{1}{\eta}=\frac{1}{c}+\frac{1}{\theta};
\]
    \item[(iv)] For $\delta$ given by

\[
2s = N\left(\frac{1}{\eta} - \frac{1}{\delta}\right),
\]
the map $ w\rightarrow \mathcal{L}^{-1}w$ is bounded from $L^\eta(\Omega)$ to $L^\delta(\Omega)$.
    \end{itemize}

Joining $(i)$, $(ii)$, $(iii)$ and $(iv)$, one easily checks that $\gamma < \delta$ and, in addition,

\begin{eqnarray*}
\Vert u\Vert_{L^\delta} &\leq & \Vert \mathcal{L}^{-1}\left[q_\varepsilon(x)\left(\mathcal{L}^{-1}(bu)\right)^{\frac{p}{2}}\right]\Vert_{L^\delta}+\Vert \mathcal{L}^{-1}f_\varepsilon\Vert_{L^\delta}\\
&\leq & C \left(\Vert q_\varepsilon\Vert_{L^{c}}\Vert u\Vert_{L^\delta}^{\frac{p}{2}} + \Vert f_\varepsilon\Vert_{L^\delta}\right).
\end{eqnarray*}
Using now the fact that $p \leq 1$, $\Vert q_\varepsilon \Vert_{L^{c}} < \varepsilon$ and $f_\varepsilon \in L^\infty(\Omega)$, we deduce that $\Vert u \Vert_{L^\delta} \leq C$ for some constant $C > 0$ independent of $u$. Proceeding inductively, we get $u \in L^\delta(\Omega)$ for all $\delta \geq 1$. Now, we use Holder's inequality to obtain
\[
\int\limits_{\Omega}\left\vert\frac{u^{q}}{\vert x\vert^{\beta}}\right\vert^{N/2s}dx\leq\left(\int\limits_{\Omega}u^{\delta}dx\right)^{qN/2s\gamma}\left(\int\limits_{\Omega}\vert x\vert^{-\beta N\delta/(2s\delta-qN)}dx\right)^{(2s\delta-qN)/2s\delta}\leq C\left(\int\limits_{\Omega}u^{\delta}dx\right)^{qN/2s\delta}
\]
if 
\[
\frac{\beta N\delta}{2s\delta-qN}<N.
\]
That is
\[
Nq<(2s-\beta)\delta.
\]
Since $u \in L^\delta(\Omega)$, for all $\delta \geq 1$, we have $\frac{u^{q}}{\vert x\vert^{\beta}}\in L^{N/2s}(\Omega)$ if $\beta<2s$. Thus, by regularity result (see Proposition 1.4 of \cite{ROS}) we have $v\in C^{\xi}(\R^N)$ for some $\xi\in(0,1)$. Analogously, we use Holder's inequality to obtain $\frac{v^{p}}{\vert x\vert^{\alpha}}\in L^{N/2s}(\Omega)$ if $\alpha<2s$. Thus, by regularity result (see Proposition 1.4 of \cite{ROS}) we have $u\in C^{\sigma}(\R^N)$ for some $\sigma\in(0,1)$.

Now if $q \leq 1$ write $b(x) = u(x)^{\frac{q}{2}}$ and if $p > 1$ write $a(x) = v(x)^{p-1}$. \ \fim

\end{document}